\documentclass[]{amsproc}

\usepackage{amsmath,amssymb,amsthm}

\newtheorem{theorem}{Theorem}
\newtheorem{lemma}{Lemma}
\newtheorem*{Bott}{Bott's theorem}

\theoremstyle{definition}
\newtheorem{example}{Example}

\theoremstyle{remark}
\newtheorem{remark}{Remark}

\begin{document}
\title[Notes on homogeneous vector bundles over complex flag manifolds]
{Notes on homogeneous vector bundles \\ over complex flag manifolds}
\author[Sergei Igonin]{Sergei Igonin}

\thanks{Partially supported by the RFBR grants 01-01-00709, 00-01-00705.}
\address{Independent University of Moscow; University of Twente, the Netherlands}
\email{igonin@mccme.ru}
\begin{abstract}
Let $P$ be a parabolic subgroup of a semisimple complex Lie group $G$
defined by a subset $\Sigma$ of simple roots of $G$,
and let $\mathbf{E}_\varphi$ be a homogeneous vector bundle over
the flag manifold $G/P$ corresponding to a linear 
representation $\varphi$ of $P$.
Using Bott's theorem, we obtain sufficient conditions
on $\varphi$ in terms of the combinatorial structure of $\Sigma$
for some cohomology groups of 
the sheaf of holomorphic sections of ${\mathbf E}_\varphi$ to be zero. 
In particular, 
we define two numbers $d(P),\,\ell(P)\in\mathbb N$ such that for any $\varphi$
obtained by natural operations from a representation of dimension
less than $d(P)$ the $q$-th cohomology group of ${\mathbf E}_\varphi$ 
is zero for $0<q<\ell(P)$. We prove also that in this case that 
the vector bundle $\mathbf{E}_\varphi$ is rigid. 
\end{abstract}
\maketitle

Let $\mathbf E$ be a holomorphic vector bundle over a connected compact
complex manifold $M$. Then there is a natural homomorphism of complex Lie
groups $\mu\colon \operatorname{Aut}{\mathbf E}\to\operatorname{Bih} M$, where
$\operatorname{Aut}{\mathbf E}$ is the automorphism group of the vector bundle and
$\operatorname{Bih} M$ is the group of all biholomorphic transformations of $M$.
The bundle $\mathbf E$ is said to be {\it homogeneous} if the action
$\mu$ of $\operatorname{Aut}\mathbf E$ on $M$ is transitive.

Assume that we have a homomorphism $\Phi\colon
G\to\operatorname{Aut}\mathbf E$ such that
the action $\mu\Phi$ of $G$ on $M$ is transitive.
Then $\mathbf E$ is homogeneous;
we say that $\mathbf E$ is {\it homogeneous with respect to $G$}.
Let $o\in M$ and let $P = G_o$ be the stabilizer of $o$ in $G$.
Then $P$ acts in a natural way on the fibre $E=E_o$, that is,
we have a holomorphic linear representation $\varphi\colon P\to\operatorname{GL}(E)$.
It is known that the bundle $\mathbf E$ is uniquely determined by the group
$G$, the subgroup $P$, and the representation $\varphi$, which can all be
arbitrary. We denote by ${\mathbf E}_{\varphi}$ the homogeneous vector
bundle over $M=G/P$ defined by a representation $\varphi$ of $P$. 

Note also that corresponding to standard tensor operations on 
representations there are similar operations on vector bundles with a 
fixed base. In particular,
\begin{equation}\label{rep&bundle}
{\mathbf E}_{\varphi^*}={\mathbf E}_{\varphi}^*,\ \
{\mathbf E}_{\varphi_1+\varphi_2}={\mathbf E}_{\varphi_1}\oplus{\mathbf E}_{\varphi_2},
\ \
{\mathbf E}_{\varphi_1\varphi_2}={\mathbf E}_{\varphi_1}\otimes{\mathbf E}_{\varphi_2}.
\end{equation}

For a Lie group $G$ denote by $G^\circ$ the identity component of $G$.
We consider the case when $M$ is a flag manifold, that is,
when $M$ is homogeneous and the stabilizers
$(\operatorname{Bih} M)^{\circ}_x\subset (\operatorname{Bih} M)^{\circ},\; x\in M$, are
parabolic subgroups (see \cite{akh,topology}).
Then the group $\operatorname{Bih} M$ and all its transitive subgroups are
semisimple Lie groups. Let ${\mathbf E}\to M$ be a homogeneous vector bundle
and let $Q=\mu((\operatorname{Aut}{\mathbf E})^{\circ})$. Then $Q$
acts on $M$ transitively
and, therefore, is semisimple. By Levi's theorem, there exists a connected
Lie subgroup $G$ of $\operatorname{Aut}{\mathbf E}$ such that $\mu\colon G\to Q$ is
a local isomorphism.
Thus ${\mathbf E}$ is homogeneous with respect to a connected semisimple
Lie group $G$ whose action on $M$ is locally effective.

We choose a maximal torus $T$ of $G$ and denote by ${\mathfrak t}$
the corresponding Cartan subalgebra
of the tangent algebra ${\mathfrak g}$ of $G$. Let
$\Delta\subset{\mathfrak t}^*$ be the root system
associated with $T$ and let $\Delta_+\subset\Delta$ be a subset of positive roots.
We denote by $\Pi$ the corresponding system of simple roots.

As usually, consider a non-degenerate $G$-invariant scalar
product in ${\mathfrak g}$ inducing a non-degenerate
scalar product $(\cdot,\cdot)$ in ${\mathfrak t}^*$ invariant under
the Weyl group $W$.
An element $\beta\in{\mathfrak t}^*$ is called a {\it weight} if 
\begin{equation}
\label{weight}
\frac{2(\beta,\alpha)}{(\alpha,\alpha)}\in\mathbb Z\,\ 
\text{ for all }\,\alpha\in\Delta.
\end{equation}
A weight $\beta$ is said to be {\it dominant}
if $(\beta,\alpha)\ge 0$ for each $\alpha\in\Delta^+$.

Let ${\mathfrak g}_{\alpha}$ be the root subspace of ${\mathfrak g}$
corresponding to $\alpha\in\Delta$. If $R$ is a Lie subgroup of $G$
normalized by $T$, then its tangent subalgebra
${\mathfrak r}\subset{\mathfrak g}$ has the following form:
$$
{\mathfrak r}=
({\mathfrak r}\cap{\mathfrak t})\oplus
\bigoplus_{\alpha\in\Delta(R)}{\mathfrak g}_{\alpha},
$$
where $\Delta(R)\subset\Delta$ is some subset called the {\it root system of} $R$.
In particular, each system of positive roots $\Delta_+$ determines two Borel
(that is, maximal solvable) subgroups $B_{\pm}$ of $G$ containing $T$
with root systems $\Delta(B_{\pm})=\pm\Delta_+$.

Let $P$ be a parabolic subgroup of $G$, that is, a subgroup containing
a Borel subgroup. We assume that $P$ contains the subgroup
$B_-$ corresponding to the system of negative roots $\Delta_- = -\Delta_+$.
It is known (see \cite{akh,topology}) that
$\Delta(P) = \Delta_-\cup [\Sigma]$ in this case, where
$[\Sigma]$ is the set of all roots of $G$ that can be expressed as a linear
combination of elements of a subset $\Sigma\subset\Pi$.
We have the semidirect decomposition
$$
P=H\ltimes N_-,\ \
{\mathfrak p}={\mathfrak h}\;+\!\!\!\!\!\!\supset{\mathfrak n}_-,
$$
where $H$ is a maximal reductive subgroup, $N_-$ is the
unipotent radical of $P$, and
${\mathfrak p},\,{\mathfrak h},\,{\mathfrak n}_-$ are the corresponding
Lie algebras. Moreover,
\begin{align*}
\Delta(H) &= [\Sigma],\\
\Delta(N_-) &= \Delta_-\setminus [\Sigma],
\end{align*}
and $\Sigma$ coincides with the system of simple roots of the group $H$
corresponding to its Borel subgroup $B_+\cap H$.
Thus the parabolic subgroup $P\subset G$ is uniquely determined
by the subset $\Sigma\subset\Pi$ of simple roots.

Note that $T$ is also a maximal torus of $P$. Each weight of a representation
$\varphi$ of $P$ is a weight in the above sense. A representation $\varphi$ is
completely reducible if and only if it is trivial on $N_-$; in this case it
is uniquely determined by the representation $\varphi|_H$ of $H$. A {\it highest weight}
of $\varphi$ is, by definition, a highest weight of $\varphi|_H$ understood in the
sense of the ordering corresponding to the Borel subgroup $B_+\cap H$.

One uses Bott's well-known theorem (see \cite{akh,bott}) to calculate
the graded cohomology space $H^*(M,\mathcal{E}_\varphi)$, where $M=G/P$ is a flag
manifold and $\mathcal{E}_\varphi$
is the sheaf of holomorphic sections of the homogeneous vector bundle
${\mathbf E}_\varphi\to M$ defined by a representation 
$\varphi\colon P\to\operatorname{GL}(E)$.

A weight $\lambda\in{\mathfrak t}^*$ is said to be {\it singular} if there exists
a root $\alpha\in\Delta_+$ such that $(\lambda,\alpha)=0$, and it is said to
be {\it regular} otherwise. We set $\gamma =\frac12\sum_{\alpha\in\Delta_+}\alpha$;
then
\begin{equation}\label{gamma}
(\gamma,\alpha)=\frac12(\alpha,\alpha)\,\ \text{ for each }\,\alpha\in\Pi.
\end{equation}
If $\lambda$ is a regular weight, then there exists a unique
$w\in W$ such that $w(\lambda)-\gamma$ is dominant. The {\it index} of the
regular weight $\lambda$ is the smallest integer $s$ such that $w$ expands into
the product of $s$ reflections $s_{\alpha},\;\alpha\in\Pi$.
This is equal to the number of $\alpha\in\Delta_+$ such that $(\lambda,\alpha)<0$.
To each weight $\lambda\in{{\mathfrak t}}^*$ such that $\lambda+\gamma$
is regular we assign a dominant weight $I(\lambda)\in{\mathfrak t}^*$
by the rule $I(\lambda)=w(\lambda+\gamma)-\gamma$.

To formulate Bott's theorem, note that in each cohomology space
$H^q(M,\mathcal{E})$ associated with a homogeneous vector bundle
${\mathbf E}\to M=G/P$ there exists a natural structure of a $G$-module.

\begin{Bott}
Let $\varphi$ be an irreducible finite-dimensional representation of
$P$ with highest weight $\Lambda$.
Then the graded cohomology space $H^*(M,\mathcal{E}_{\varphi})$
is determined by $\Lambda$ in the following way:
\begin{itemize}
\item if $\Lambda +\gamma$ is singular then $H^*(M,\mathcal{E}_{\varphi})=0$;
\item if $\Lambda +\gamma$ is regular and its index is $p$,
then $H^q(M,\mathcal{E}_{\varphi})=0$ for all $q\ne p$, while 
$H^p(M,\mathcal{E}_{\varphi})$
is an irreducible $G$-module with highest weight $I(\Lambda)$.
\end{itemize}
\end{Bott}
Using this theorem, we will obtain some sufficient conditions on
representations $\varphi$ of $P$ for a cohomology group
$H^q(M,\mathcal{E}_{\varphi})$ to be zero.

Bott's theorem is not applicable directly when the representation
is not completely reducible. But due to the following lemma 
to prove that cohomology is zero
it is enough to consider completely reducible representations.
For an arbitrary holomorphic representation $\varphi:P\to\operatorname{GL}(E)$ we 
construct a completely reducible representation $\varphi^s$ as follows. 
There is a filtration ({\it Jordan-H\"older tower})
$$
0=E_0\subset E_1\subset\ldots\subset E_m=E
$$ 
of $E$ 
by invariant subspaces of $\varphi(P)$ such that for all $i=1,\ldots,m$ the 
induced representation $\bar\varphi_i$ in $E_i/E_{i-1}$ is irreducible.
We set 
\begin{equation}
  \label{ph_s}
  \varphi^s=\bar\varphi_1+\dots+\bar\varphi_m.
\end{equation}
It is well-known that $\varphi^s$ does not depend on a Jordan-H\"older
tower of $E$. 
\begin{lemma}\label{reducible}
If $H^q(M,\mathcal{E}_{\varphi^s})=0$ for some $q\ge 0$ then 
$H^q(M,\mathcal{E}_{\varphi})=0$.
\end{lemma}
\begin{proof}
From $H^q(M,\mathcal{E}_{\varphi^s})=0$ and (\ref{ph_s}) we have
\begin{equation}\label{factor}
H^q(M,\mathcal{E}_{\bar\varphi_i})=0\,\ \text{ for all }\,i=1,\ldots,m.
\end{equation}
Denote by $\varphi_i$ the restriction of $\varphi$ onto the invariant
subspace $E_i$ and consider the short exact sequences of sheaves
$$
0\to\mathcal{E}_{\varphi_{i-1}}\to\mathcal{E}_{\varphi_i}\to
\mathcal{E}_{\bar\varphi_i}\to 0,\,\ 1\le i\le m,
$$
corresponding to the exact sequences
$$
0\to E_{i-1}\to E_i\to E_i/E_{i-1}\to 0
$$
of $P$-modules.
Using the exact cohomology sequence, by induction on $i$ from
(\ref{factor}) we get $H^q(M,\mathcal{E}_{\varphi_i})=0$ for all $i=1,\ldots,m$.
\end{proof}%\DeclareMathOperator{\im}{im}

Each $\xi\in\mathfrak{t}^*$ is a linear combination of simple roots
$\alpha\in\Pi$. We denote by $c_\alpha(\xi)$ the corresponding 
coefficients, that is, $\xi=\sum_{\alpha\in\Pi}c_\alpha(\xi)\cdot\alpha$,
and set
$$
\mathrm{C}(\xi)=\{\alpha\in\Pi: c_\alpha(\xi)\neq 0\}\subset\Pi.
$$
The value $\sqrt{(\xi,\xi)}$ is called the {\it length} of $\xi$ 
and denoted by $|\xi|$.
For $\delta^1,\,\delta^2\in\Delta_+$ we write $\delta^1\le\delta^2$ if 
$c_\alpha(\delta^1)\le c_\alpha(\delta^2)$ for all $\alpha\in\Pi$. 
We need two lemmas on semisimple Lie algebras.
\begin{lemma}
  \label{root}
If $\delta^1,\,\delta^2\in\Delta_+$ and $\delta^1\le\delta^2$ then there exists
a sequence of positive roots $\delta_0,\,\delta_1,\ldots,\,\delta_m$
such that $\delta_0=\delta^1$, $\delta_m=\delta^2$, and 
$\delta_i-\delta_{i-1}\in\mathrm{C}(\delta^2-\delta^1)$ for all $i=1,\ldots,m$.
\end{lemma}
\begin{proof}
By induction on $|\mathrm{C}(\delta^2-\delta^1)|$, we must prove that
if $\delta^1\neq\delta^2$ then there is $\alpha\in\mathrm{C}(\delta^2-\delta^1)$ 
such that
\begin{equation}
  \label{or}
  \delta^2-\alpha\in\Delta_+\, \ \text{ or }\ \, \delta^1+\alpha\in\Delta_+.
\end{equation}
We have $(\delta^2-\delta^1,\delta^2-\delta^1)>0$, hence 
$(\delta^2-\delta^1,\alpha)>0$ for some
$\alpha\in\mathrm{C}(\delta^2-\delta^1)$. 
Then $(\delta^2,\alpha)>0$ or $(\delta^1,\alpha)<0$, which implies (\ref{or}). 
\end{proof}
\begin{lemma}\label{dimless}
Let $\mathfrak{g}_1,\ldots,\mathfrak{g}_n$
be simple Lie algebras. Denote by $d_i$ 
the minimal dimension of a nontrivial representation of $\mathfrak{g}_i$.
Let $\psi$ be a representation of the semisimple algebra $\oplus_i\mathfrak{g}_i$. 
If $\dim\psi<\sum_i d_i$ then $\psi$ is trivial on some $\mathfrak{g}_i$.
\end{lemma}
\begin{proof}
Each irreducible component of $\psi$ is of the form
\begin{equation}\label{otimes}
\psi_{i_1}\otimes\psi_{i_2}\otimes\dots\otimes\psi_{i_k},
\end{equation}
where $\psi_{i_j}$ is a nontrivial irreducible representation of
$\mathfrak{g}_{i_j}$. The dimension of (\ref{otimes}) is not less than
$$
d_{i_1} d_{i_2}\dots d_{i_k}\ge
d_{i_1}+d_{i_2}+\dots+d_{i_k}.
$$
Therefore, if $\psi$ is nontrivial on each $\mathfrak{g}_i$ then
$\dim\psi\ge\sum_i d_i$.
\end{proof}

Let $\mathrm{A}\subset\Pi\setminus\Sigma$ and $\mathrm{B}\subset\Sigma$.
In the sequel we consider representations $\varphi$ of $P$ with
highest weights $\Lambda$ satisfying
\begin{equation}
  \label{AB}
   (\Lambda,\alpha)<0,\ (\Lambda,\beta)=0\,\ \text{ for all }\,
  \alpha\in\mathrm{A},\,\beta\in\mathrm{B}.
\end{equation}
We define a class of positive roots $\delta$ of $\mathfrak{g}$
({\it significant $(\mathrm{A},\mathrm{B})$-roots}) such that
if $\Lambda+\gamma$ is regular then $(\Lambda+\gamma,\delta)<0$.
By Bott's theorem and Lemma \ref{reducible}, this implies
that $H^q(M,\mathcal{E}_{\varphi})=0$ for $0\le q<\ell(\mathrm{A},\mathrm{B})$, 
where $\ell(\mathrm{A},\mathrm{B})$ is the number of significant 
$(\mathrm{A},\mathrm{B})$-roots.

Let us introduce the required concepts. A positive root
$$
\delta=\sum_{\alpha\in\mathrm{A}}c_\alpha(\delta)\cdot\alpha+
\sum_{\beta\in\mathrm{B}}c_\beta(\delta)\cdot\beta,\,\
c_\alpha(\delta),\,c_\beta(\delta)\ge 0,
$$
of $\mathfrak{g}$ with 
\begin{equation}
  \label{c_al}
  c_{\alpha_0}(\delta)>0\,\ \text{ for some }\,\alpha_0\in\mathrm{A}
\end{equation}
will be called an {\it $(\mathrm{A},\mathrm{B})$-root}. Denote by 
$\Delta_{\mathrm{A},\mathrm{B}}$
the set of $(\mathrm{A},\mathrm{B})$-roots.
We say that $\delta\in\Delta_{\mathrm{A},\mathrm{B}}$ is {\it significant} 
if there exists a root $\sigma_\delta\in\Delta_{\mathrm{A},\mathrm{B}}$ such that
\begin{equation}
  \label{sign1}
  \sigma_\delta\le\delta,\ \ 
  \sum_{\alpha\in\mathrm{A}}c_\alpha(\sigma_\delta)\cdot(\alpha,\alpha)\ge
  \sum_{\beta\in\mathrm{B}}c_\beta(\sigma_\delta)\cdot(\beta,\beta),   
\end{equation}
and $\mathrm{C}(\delta-\sigma_\delta)$ consists of simple roots of the same length 
not greater than $|\sigma_\delta|$.

\begin{lemma}
\label{significant}
If a weight $\Lambda$ satisfies (\ref{AB}) and $\Lambda+\gamma$ is regular
then for each significant $\delta\in\Delta_{\mathrm{A},\mathrm{B}}$ one has
$(\Lambda+\gamma,\delta)<0$.
\end{lemma}
\begin{proof}
By (\ref{AB}) and (\ref{gamma}), we have
\begin{equation}
\label{long}
(\Lambda+\gamma,\sigma_\delta)=
\sum_{\alpha\in\mathrm{A}}(c_\alpha(\sigma_\delta)\cdot(\Lambda,\alpha)
+\frac12 c_\alpha(\sigma_\delta)\cdot(\alpha,\alpha))+
\sum_{\beta\in\mathrm{B}}\frac12 c_\beta(\sigma_\delta)\cdot(\beta,\beta).
\end{equation}
Note that 
\begin{equation}
\label{-1}
\frac{2(\Lambda,\alpha)}{(\alpha,\alpha)}\neq -1\,\ \text{ for each }\,\alpha\in\Pi
\end{equation}
because otherwise $(\Lambda+\gamma)(\alpha)=0$ and $\Lambda+\gamma$ is singular.
Combining (\ref{weight}), (\ref{AB}), and (\ref{-1}), we obtain
\begin{equation}
\label{-2}
\frac{2(\Lambda,\alpha)}{(\alpha,\alpha)}\le -2\,\ 
\text{ for each }\,\alpha\in\mathrm{A}.
\end{equation}
From (\ref{long}), (\ref{-2}), and (\ref{sign1}) one gets
\begin{equation}
\label{0<0}
(\Lambda+\gamma,\sigma_\delta)\le
\sum_{\alpha\in\mathrm{A}}-\frac12c_\alpha(\sigma_\delta)\cdot(\alpha,\alpha)+
\sum_{\beta\in\mathrm{B}}\frac12 c_\beta(\sigma_\delta)\cdot(\beta,\beta)\le 0.
\end{equation}
By Lemma \ref{root}, there exist 
$\delta_0,\,\delta_1,\ldots,\,\delta_m\in\Delta_{\mathrm{A},\mathrm{B}}$ such 
that $\delta_0=\sigma_\delta$, $\delta_m=\delta$, and 
\begin{equation}
\label{inC}
\delta_i-\delta_{i-1}\in\mathrm{C}(\delta-\sigma_\delta)\subset\mathrm{A}\cup\mathrm{B}.
\end{equation}
Denote 
$$
c=\frac12(\delta_i-\delta_{i-1},\delta_i-\delta_{i-1}).
$$ 
By the definition of significant roots, 
$c$ does not depend on $i=1,\ldots,m$ and 
\begin{equation}
  \label{cle}
  2c\le(\sigma_\delta,\sigma_\delta).
\end{equation}
For any roots $\alpha,\,\beta$ if $|\alpha|\ge|\beta|$ then
$\frac{(\alpha,\alpha)}{(\beta,\beta)}\in\mathbb N$. Applying this to
$\delta_0=\sigma_\delta$ and $\delta_i-\delta_{i-1}$, from (\ref{cle}) we get
$\frac{1}{2c}(\delta_0,\delta_0)\in\mathbb N$
and, therefore, from (\ref{weight}) one obtains
\begin{equation}
  \label{delta0}
  \frac1c(\Lambda+\gamma,\delta_0)=
\frac{2(\Lambda+\gamma,\delta_0)}{(\delta_0,\delta_0)}
\cdot\frac{(\delta_0,\delta_0)}{2c}\in\mathbb Z.
\end{equation}
By (\ref{weight}) and (\ref{gamma}),
\begin{equation}
\label{inZ}
\frac1c(\Lambda+\gamma,\delta_i-\delta_{i-1})=
\frac1c(\Lambda,\delta_i-\delta_{i-1})+1\in\mathbb Z.
\end{equation}
Moreover, from (\ref{inC}) and (\ref{AB}) we have 
$\frac1c(\Lambda,\delta_i-\delta_{i-1})\le 0$, hence
\begin{equation}
  \label{le1}
  \frac1c(\Lambda+\gamma,\delta_i-\delta_{i-1})\le 1.
\end{equation}
By induction on $i$, combining (\ref{0<0}), 
(\ref{delta0}), (\ref{inZ}), and (\ref{le1}), we get
$\frac1c(\Lambda+\gamma,\delta_i)\in\mathbb Z$ and 
\begin{equation}
\label{i<0}
(\Lambda+\gamma,\delta_i)<0 
\,\ \text{ for all }\,i=0,1,\ldots,m
\end{equation}
because otherwise $(\Lambda+\gamma,\delta_i)=0$ for some 
positive root $\delta_i$ and $\Lambda+\gamma$ is singular.
In particular, $(\Lambda+\gamma,\delta)<0$. 
\end{proof}
\begin{example}
  \label{significant1}
We prove three sufficient conditions for $\delta\in\Delta_{\mathrm{A},\mathrm{B}}$ 
to be significant by giving corresponding $\sigma_\delta$. 
\begin{enumerate}
\item
{\it If $\mathrm{C}(\delta)$ contains simple roots of the same length},
we can, by (\ref{c_al}), take $\sigma_\delta\in\mathrm{C}(\delta)\cap\mathrm{A}$.
\item
{\it If there is $\alpha_0\in\mathrm{C}(\delta)\cap\mathrm{A}$ longer
than the other roots from $\mathrm{C}(\delta)\setminus\{\alpha_0$\}}, then 
$\mathrm{C}(\delta)\setminus\{\alpha_0\}$ contains roots of the same length,
$c_{\alpha_0}(\delta)=1$, hence 
$\mathrm{C}(\delta-\alpha_0)=\mathrm{C}(\delta)\setminus\{\alpha_0\}$ and one can put
$\sigma_\delta=\alpha_0$.
\item
{\it If there exists $\alpha_0\in\mathrm{C}(\delta)\cap\mathrm{A}$ shorter
than the other roots from $\mathrm{C}(\delta)\setminus\{\alpha_0\}$ and
$c_{\alpha_0}(\delta)=-\frac{2(\beta,\alpha_0)}{(\alpha_0,\alpha_0)}$ for some
$\beta\in\mathrm{C}(\delta)\setminus\{\alpha_0$\}},
then $\mathrm{C}(\delta)\setminus\{\alpha_0\}$ consists of roots of 
length $|\beta|=|\beta-\frac{2(\beta,\alpha_0)}{(\alpha_0,\alpha_0)}\cdot\alpha_0|$ and
we set 
$$
\sigma_\delta=\beta-\frac{2(\beta,\alpha_0)}{(\alpha_0,\alpha_0)}\cdot\alpha_0.
$$
\end{enumerate}
\end{example}

Let $\ell(\mathrm{A},\mathrm{B})$ be the number of significant 
$(\mathrm{A},\mathrm{B})$-roots.
This number can be computed in terms of some graphs,
if we regard $\mathrm{A}$ and $\mathrm{B}$ as subgraphs of the Dynkin diagram
$\mathrm{D}$ of $\mathfrak g$.
\begin{example}\label{length}
Denote by $\mathrm{U}$ the union of the connected components $\mathrm{C}$
of $\mathrm{A}\cup\mathrm{B}$ satisfying
$\mathrm{C}\cap\mathrm{A}=\emptyset$ and set 
$\mathrm{B}'=\mathrm{B}\setminus\mathrm{U}$. 
Clearly, $(\mathrm{A},\mathrm{B})$-roots coincide with $(\mathrm{A},\mathrm{B}')$-roots and  
$\ell(\mathrm{A},\mathrm{B})=\ell(\mathrm{A},\mathrm{B}')$ because for
each $(\mathrm{A},\mathrm{B})$-root $\delta$ 
the subgraph $\mathrm{C}(\delta)\subset\mathrm{A}\cup\mathrm{B}$ is connected and, 
by (\ref{c_al}), $\mathrm{C}(\delta)\cap\mathrm{A}\neq\emptyset$.

Suppose that all the edges of $\mathrm{A}\cup\mathrm{B}'$ are simple, that is, 
each connected component of $\mathrm{A}\cup\mathrm{B}'$ contains simple roots of 
the same length. Then, by Example \ref{significant1}, each $(\mathrm{A},\mathrm{B}')$-root 
is significant, 
hence $\ell(\mathrm{A},\mathrm{B})$ equals the number of $(\mathrm{A},\mathrm{B}')$-roots.
The latter is equal to the number of positive roots of 
$\mathfrak{g}_{\mathrm{A}\cup\mathrm{B}'}$ 
minus the number of positive roots of $\mathfrak{g}_{\mathrm{B}'}$,
where $\mathfrak{g}_{\mathrm{A}\cup\mathrm{B}'}$ and
$\mathfrak{g}_{\mathrm{B}'}$ 
are semisimple Lie algebras with
Dynkin diagrams $\mathrm{A}\cup\mathrm{B}'$ and $\mathrm{B}'$ respectively. 
\end{example}

In what follows, we will denote by the same letter 
a representation $\varphi\colon P\to\operatorname{GL}(E)$ and its differential
$\varphi\colon\mathfrak{p}\to\mathfrak{gl}(E)$.
\begin{theorem}\label{main}
Consider an irreducible representation $\varphi$ of $P$ with
highest weight $\Lambda\in{\mathfrak t}^*$. Let
$\mathrm{A}$ and $\mathrm{B}$ be subsets of $\Pi\setminus\Sigma$ and $\Sigma$ respectively
such that (\ref{AB}) holds.
Then we have $H^q(M,\mathcal{E}_{\varphi})=0$ for $0\le q<\ell(\mathrm{A},\mathrm{B})$.
\end{theorem}
\begin{proof}
If $\Lambda+\gamma$ is singular then, by Bott's theorem, 
$H^*(M,\mathcal{E}_{\varphi})=0$; suppose that $\Lambda +\gamma$ is regular. 
By Lemma \ref{significant}, $(\Lambda +\gamma,\delta)<0$ for
each significant $\delta\in\Delta_{\mathrm{A},\mathrm{B}}$.
Therefore, the index of $\Lambda+\gamma$ is not less than $\ell(\mathrm{A},\mathrm{B})$. 
Applying Bott's theorem, one completes the proof.
\end{proof}

We say that a representation $\varphi\colon P\to\operatorname{GL}(E)$ {\it is obtained
by natural operations from a representation
$\tilde\varphi\colon P\to\operatorname{GL}(\tilde E)$} if
there is a homomorphism $\pi\colon\operatorname{GL}(\tilde
E)\to\operatorname{GL}(E)$ such that
\begin{equation}\label{natural}
\varphi=\pi\circ\tilde\varphi.
\end{equation}
\begin{example}\label{tensor}
A representation $\varphi$ obtained from $\tilde\varphi$ by tensor operations 
(\ref{rep&bundle}), for instance $\varphi=\tilde\varphi{\tilde\varphi}^*$, 
satisfies (\ref{natural}).
\end{example}
Denote by $H'$ and ${\mathfrak h}'$ the semisimple commutator subgroup of $H$
and its tangent subalgebra respectively. 
Let us regard $\Sigma$ as a subdiagram of the Dynkin diagram $\mathrm{D}$ of
$\mathfrak g$. Then
${\mathfrak h}'$ is a semisimple Lie algebra with Dynkin diagram $\Sigma$.
As usually, each connected component $\mathrm{C}$ of $\Sigma$ determines a simple ideal
${\mathfrak h}_\mathrm{C}$ of ${\mathfrak h}'$.
If a representation $\tilde\varphi$ of $P$ is trivial on some 
${\mathfrak h}_\mathrm{C}$ then, by (\ref{natural}), any $\varphi$ obtained from 
$\tilde\varphi$ by natural operations is also trivial on ${\mathfrak h}_\mathrm{C}$.
This observation, Lemmas \ref{reducible},\ref{dimless}, 
and Theorem \ref{main} allow to prove that if $\varphi$ is obtained by natural
operations from a representation of relatively small dimension then
some cohomology groups $H^q(M,\mathcal{E}_\varphi)$ are zero. 

Indeed, Lemma \ref{dimless} guarantees that $\tilde\varphi$ and, therefore, 
$\varphi$ are trivial on some ideals
of $\mathfrak{h}'$ if the dimension of $\tilde\varphi$ is sufficiently small. 
The less is the dimension of $\tilde\varphi$ the more ideals $\varphi$ is
trivial on. By the construction (\ref{ph_s}) of the corresponding 
completely reducible representation $\varphi^s$, $\varphi(\mathfrak{h}_\mathrm{C})=0$ implies
$\varphi^s(\mathfrak{h}_\mathrm{C})=0$. Let $\mathrm{B}$ be the union of those connected 
components $\mathrm{C}$ of $\Sigma$ for which one has $\varphi^s(\mathfrak{h}_\mathrm{C})=0$.

Consider an arbitrary irreducible component $\varphi'$ of $\varphi^s$ with
highest weight $\Lambda$. Clearly,
$(\Lambda,\beta)=0$ for each $\beta\in\mathrm{B}$.
We set 
$$
\mathrm{A}(\varphi')=\{\alpha\in\Pi\setminus\Sigma:(\Lambda,\alpha)<0\}.
$$
If $\mathrm{A}(\varphi')$ is empty then $\Lambda$ is dominant and, 
by Bott's theorem, $H^q(M,\mathcal{E}_{\varphi'})=0$ for all $q>0$.
Otherwise, according to Theorem \ref{main}, 
we get $H^q(M,\mathcal{E}_{\varphi'})=0$ for $0\le q<\ell(\mathrm{A}(\varphi'),\mathrm{B})$. 
Thus, by Lemma \ref{reducible},
\begin{equation}
  \label{h=h}
  H^q(M,\mathcal{E}_{\varphi^s})=H^q(M,\mathcal{E}_{\varphi})=0\,\ \text{ for all }\,
  0<q<\min_{\varphi'}\ell(\mathrm{A}(\varphi'),\mathrm{B}),
\end{equation}
where $\varphi'$ runs through the irreducible components of $\varphi^s$ 
with $\mathrm{A}(\varphi')\neq\emptyset$. 

The following theorem is an example of such a result. 
To formulate it we need some notations.
We say that a subset $\mathrm{B}\subset\Sigma$ is {\it adjacent}
to a vertex $\alpha\in\Pi\setminus\Sigma$ of $\mathrm{D}$ if there is an edge
of $\mathrm{D}$ connecting $\alpha$ and $\mathrm{B}$. Equivalently, there is a unique
simple root $\beta\in\mathrm{B}$ such that $(\alpha,\beta)\neq 0$.
For a simple root $\alpha\in\Pi\setminus\Sigma$ consider the connected
components $\mathrm{C}_1,\ldots,\mathrm{C}_n$ of
$\Sigma\subset\mathrm{D}$ adjacent to $\alpha$ 
and the corresponding simple ideals
${\mathfrak h}_{\mathrm{C}_1},\ldots,{\mathfrak h}_{\mathrm{C}_n}$
of ${\mathfrak h}'$. Denote by $d_i$ the minimal dimension
of a nontrivial representation of ${\mathfrak h}_{\mathrm{C}_i}$
and set
\begin{align*}
d(\alpha)&=d_1+\dots+d_n,\\
\ell(\alpha)&=\min_{i=1,\dots,n}\ell(\{\alpha\},\mathrm{C}_i).
\end{align*}
If there are no nonempty connected components
adjacent to $\alpha$, we put $d(\alpha)=\ell(\alpha)=1$. Set also
\begin{align*}
d(P)&=\min_{\alpha\in\Pi\setminus\Sigma}d(\alpha),\\
\ell(P)&=\min_{\alpha\in\Pi\setminus\Sigma}\ell(\alpha).
\end{align*}
\begin{theorem}\label{H1}
Suppose that a representation $\varphi$ of $P$ is obtained by natural
operations from a representation $\tilde\varphi$ of dimension less then $d(P)$.
Then we have 
\begin{equation}\label{H1=0}
H^q(M,\mathcal{E}_{\varphi})=0\,\ \text{ for }\,0<q<\ell(P).
\end{equation}
\end{theorem}
\begin{proof}
If $d(P)=1$, the statement is trivial. 
If $d(P)>1$ then from Lemma \ref{dimless} it follows
that for each $\alpha\in\Pi\setminus\Sigma$ there is a nonempty
connected component $\mathrm{C}_\alpha$ of $\Sigma$ adjacent to $\alpha$ such that
$\tilde\varphi$ and, by (\ref{natural}), $\varphi$ are trivial on
${\mathfrak h}_{\mathrm{C}_\alpha}$. By the definition of $\ell(P)$,
we have $\ell(P)\le\ell(\{\alpha\},\mathrm{C}_\alpha)$ and
from (\ref{h=h}) one obtains (\ref{H1=0}).
\end{proof}
\begin{remark}
  \label{ell(P)}
We have  $d(P)>1$ if and only if for each 
$\alpha\in\Pi\setminus\Sigma$ there exists a nonempty
connected component $\mathrm{C}\subset\Sigma$ adjacent to $\alpha$.
According to the following lemma, this is also equivalent to
$\ell(P)>1$. In this case the result of Theorem \ref{H1} is
nontrivial.
\end{remark}
\begin{lemma}
  \label{a+b}
If $\mathrm{B}$ is a nonempty subset of $\Sigma$ adjacent to 
$\alpha\in\Pi\setminus\Sigma$ then 
$$
\ell(\{\alpha\},\mathrm{B})\ge 2.
$$
\end{lemma}
\begin{proof}
The root $\alpha$ is clearly a significant $(\alpha,\mathrm{B})$-root.
By assumption, there is $\beta\in\mathrm{B}$ such that 
$(\alpha,\beta)\neq 0$. Consider the root $\delta$ of the form 
$$
\delta=a\cdot\alpha+b\cdot\beta,\ \ a,\,b>0,
$$
such that $\delta+\alpha$ and $\delta+\beta$ are not roots of $\mathfrak{g}$.
According to Example \ref{significant1}, $\delta$ is also 
significant. Thus we have at least two distinct significant roots.
\end{proof}
\begin{example}
Consider the Grassmanian ${\mathrm{Gr}}_{n,k}$ of $k$-dimensional subspaces of
${\mathbb C}^n$. It is a flag manifold of the group $\operatorname{SL}_n({\mathbb C})$.
The stabilizer of a point $o\in{\mathrm{Gr}}_{n,k}$ is a parabolic subgroup with
$\Sigma=\Pi\setminus\{\alpha\}$, where $\alpha$ is the $k$-th simple root in the
Dynkin diagram
$$
\begin{picture}(200,37)
\put(40,17){\circle{6}}
\put(43,17){\line(1,0){20}}
\put(66,17){\circle{6}}
\put(69,17){\line(1,0){10}}
\put(99,17){\line(1,0){10}}
\put(112,17){\circle*{6}}
\put(105,2){$k$\text{-th root}}
\put(115,17){\line(1,0){10}}
\multiput(84,17)(5,0){3}{\circle*{2}}
\put(75,27){$\mathrm{C}_1$}
\multiput(130,17)(5,0){3}{\circle*{2}}
\put(139,27){$\mathrm{C}_2$}
\put(145,17){\line(1,0){10}}
\put(158,17){\circle{6}}
\put(161,17){\line(1,0){20}}
\put(184,17){\circle{6}}
\end{picture}
$$
of $\operatorname{SL}_n({\mathbb C})$. The subdiagram $\Sigma$
consists of two components $\mathrm{C}_1$ and
$\mathrm{C}_2$ adjacent to $\alpha$,
one of which is empty if ${\mathrm{Gr}}_{n,k}\simeq\mathbb{CP}^{n-1}$,
that is, if $k=1$ or $k=n\!-\!1$. Hence $d(P)=n$ if $1<k<n\!-\!1$, and
$d(P)=n\!-\!1$ if $k=1$ or $k=n\!-\!1$.
According to Example \ref{length}, we have $\ell(\{\alpha\},\mathrm{C}_1)=k$,
$\ell(\{\alpha\},\mathrm{C}_2)=n\!-\!k$. Thus $\ell(P)=\min\{k,n\!-\!k\}$
if $1<k<n\!-\!1$, and $\ell(P)=n\!-\!1$ if $k=1$ or $k=n\!-\!1$.
This computation gives the following result.
\begin{theorem}
For any vector bundle $\mathbf E\to{\mathrm{Gr}}_{n,k},\ 1<k<n\!-\!1$, obtained
by tensor operations from a homogeneous vector bundle 
$\tilde{\mathbf E}\to{\mathrm{Gr}}_{n,k}$ of rank less than $n$ one has 
$H^q({\mathrm{Gr}}_{n,k},\mathcal{E})=0$ for all $0<q<\min\{k,n\!-\!k\}$. 
\end{theorem}
\begin{proof}
It is known (see \cite{topology}) that for $1<k<n\!-\!1$ 
any transitive Lie subgroup of $\operatorname{Bih}{\mathrm{Gr}}_{n,k}$ coincides with
$\operatorname{Bih}{\mathrm{Gr}}_{n,k}=\operatorname{PSL}_n({\mathbb C})$. 
Therefore, each homogeneous
vector bundle over ${\mathrm{Gr}}_{n,k}$ is homogeneous with respect to the
simply connected group $\operatorname{SL}_n({\mathbb C})$.
Hence $\tilde{\mathbf E}={\mathbf E}_{\tilde\varphi}$ for some representation 
$\tilde\varphi$ of $P\subset\operatorname{SL}_n({\mathbb C})$. By (\ref{rep&bundle}), 
we have $\mathbf E={\mathbf E}_\varphi$, where $\varphi$ is
obtained from $\tilde\varphi$ by tensor operations.
According to Example \ref{tensor}, Theorem \ref{H1}, and 
the above computation of $d(P)$ and $\ell(P)$, we complete the proof.   
\end{proof}
\end{example}
Cohomology groups of sheaves considered here appear, in particular,
in deformation theory \cite{deform}. Theorems \ref{main} and \ref{H1} 
allow to prove {\it rigidity}, that is, absence of nontrivial {\it
deformations} \cite{deform}, of some homogeneous vector bundles 
and {\it supermanifolds} \cite{manin,super}. 
A result on rigidity of vector bundles is given by the following 
theorem, while homogeneous supermanifolds are studied in \cite{mine}.
\begin{theorem}
For any representation $\varphi$ of $P$ obtained by natural operations from
a representation $\tilde\varphi$ of dimension less than $d(P)$ the bundle
$\mathbf{E}_\varphi$ is rigid.
\end{theorem}
\begin{proof}
By (\ref{rep&bundle}), we have 
$\mathbf{E}_\varphi\otimes\mathbf{E}_\varphi^*=\mathbf{E}_{\varphi\varphi^*}$. 
According to Example \ref{tensor}, the representation $\varphi\varphi^*$ is
also obtained from $\tilde\varphi$ by natural operations. Then from
Theorem \ref{H1} and Remark \ref{ell(P)} one gets 
$H^1(M,\mathcal{E}_\varphi\otimes\mathcal{E}_\varphi^*)=0$, 
which implies that $\mathbf{E}_\varphi$ 
is rigid (see \cite{deform1,deform}).  
\end{proof}

\section*{Acknowledgements}
I would like to express deep gratitude to 
Prof.~A.~L.~Onishchik for his attention to this work and many 
useful remarks. It is a pleasure to thank the
organizers of the European Summer School 2001 on Asymptotic Combinatorics 
with Application to Mathematical Physics in St.~Petersburg for the
invitation and providing a stimulating environment.
I am also grateful to the Erwin Schr\"odinger Institute in Vienna for 
hospitality in October-November 2000 where a part of this research was done.

%%%%%%%%%%%%%%%%%%%%%%%%%%%%%%%%%%%%%%%%%%%%%%%%%%

\end{document}